\IfFileExists{proc-l.cls}{\documentclass{proc-l}}{\documentclass[10pt]{amsart}}
\usepackage[T1]{fontenc}
\usepackage{lmodern,microtype}
\usepackage{amsmath,amssymb,mathtools}
\usepackage{hyperref}
\hypersetup{hidelinks,pdfauthor={Nguyen Quang Dieu, Nguyen Thi Minh Thao}}
\numberwithin{equation}{section}
\newtheorem{mainthm}{Theorem}

\newtheorem{theorem}{Theorem}[section]
\newtheorem{lemma}[theorem]{Lemma}
\newtheorem{proposition}[theorem]{Proposition}
\newtheorem{corollary}[theorem]{Corollary}
\theoremstyle{remark}
\newtheorem{remark}[theorem]{Remark}
\newcommand{\PSH}{\operatorname{PSH}}
\newcommand{\MPSH}{\operatorname{MPSH}}
\newcommand{\supp}{\operatorname{supp}}
\newcommand{\ddc}{dd^c}
\newcommand{\CC}{\mathbb C}
\newcommand{\RR}{\mathbb R}
\newcommand{\ZZ}{\mathbb Z}
\newcommand{\PP}{\mathbb P}
\newcommand{\ii}{\mathrm i}
\newcommand{\loc}{\mathrm{loc}}
\newcommand{\cpt}{\Subset}
\title[Maximal plurisubharmonic functions]{Maximal plurisubharmonic functions on complex manifolds}
\author[Nguyen Quang Dieu]{Nguyen Quang Dieu}
\address{Department of Mathematics, School of Mathematics and Computer Science, Hanoi National University of Education, Hanoi, Vietnam}
\email{ngquang.dieu@hnue.edu.vn}
\author[Nguyen Thi Minh Thao]{Nguyen Thi Minh Thao}
\address{Vietnam Military Medical Academy, Hanoi, Vietnam}
\email{nguyenminhthao31502@gmail.com}
\begin{document}
\begin{abstract}
We study local and global maximality of plurisubharmonic functions on complex
manifolds. When a locally bounded, locally maximal function fails to be globally
maximal, the maximum set of the upper semicontinuous regularization of the
comparison difference is a compact pseudoconcave set. This set has no
neighborhood carrying a bounded psh function with almost everywhere positive
Monge--Amp\`ere density, and it contains a minimal pseudoconcave set on which every
positive $dd^c$-closed current has full support. Domains in simple complex tori
provide complete K\"ahler counterexamples without positive-dimensional compact
analytic subvarieties; on them, dense Levi leaves reduce comparison to convexity,
and Fourier analysis classifies the currents in dimension two. We also show that
every psh function is globally maximal on the total space of a holomorphic
submersion whose fibers are noncompact and carry only constant psh functions,
and we treat divisorial poles.
\end{abstract}
\maketitle

\section{Introduction}
Let $X$ be a connected noncompact complex manifold of dimension
$n\ge2$. A function $u\in\PSH(X)$ is \emph{globally maximal} if,
for every $v\in\PSH(X)$ and every compact set $K\subset X$,
\begin{equation}\label{eq:defmax}
 v\le u\quad\text{on }X\setminus K
 \quad\text{implies}\quad v\le u\quad\text{on }X.
\end{equation}
It is \emph{locally maximal} when the same comparison property
holds in a neighborhood of every point. Here psh functions are
not identically $-\infty$, $\MPSH(X)$ denotes the globally maximal
ones, and $d^c=\frac{\ii}{2\pi}(\bar\partial-\partial)$, so that
$\ddc=\frac{\ii}{\pi}\partial\bar\partial$.

Whether local maximality implies global maximality for \emph{arbitrary}
psh functions on domains in $\CC^n$ remains an open problem.
Originally posed by B\l ocki, it appears as Question~3 in the survey
of Dinew, Guedj and Zeriahi \cite{DinewGuedjZeriahi}.
For locally bounded functions, Bedford--Taylor theory characterizes
maximality by $(\ddc u)^n=0$ and hence answers the question
affirmatively \cite{BT82}. The same holds in the domain of
definition of the Monge--Amp\`ere operator \cite{Blocki},
and recent work treats certain locally analytic singularities
\cite{DoHiep2026}. On complex manifolds, a new obstruction appears: even
smooth functions may be locally maximal without being globally maximal.

Such examples are easy to produce. If $Z$ is
a positive-dimensional compact complex manifold, the function
$u(x,z)=|z|^2$ on $Z\times\mathbb D$ is locally maximal, but
$\max\{u,\frac14\}$ is a compactly supported competitor above $u$.
Thus compact complex directions may obstruct global comparison
without contributing to the top Monge--Amp\`ere measure. Do and Hiep
\cite{DoHiepComp} recently proved comparison principles for local
potentials on complex manifolds carrying a global psh function whose
Monge--Amp\`ere measure has a positive density almost everywhere, and
they observed the same phenomenon on $\PP^1\times\CC$. This raises two
questions: what geometry must a compact comparison defect contain,
and are compact analytic subvarieties essential to the phenomenon?
Our first main result answers the first question at two levels. We call a
compact set $K\subset X$ \emph{MA-positive} if it has a neighborhood
carrying a bounded psh function $\rho$ with
$(\ddc\rho)^n\ge f\,dV$ for some function $f>0$ almost everywhere,
where $dV$ is a smooth volume form. Recall also that a closed set is
$1$-pseudoconcave if no smooth strictly psh function has a strict
local maximum on it; see Section~\ref{sec:obstruction} for references.

\begin{mainthm}\label{thmA}
Let $u\in\PSH(X)\cap L^\infty_{\loc}(X)$ be locally maximal. Then the
following assertions hold true.
\begin{enumerate}
\item[(i)] If every compact subset of $X$ is MA-positive, then $u$ is
globally maximal. This is the case when $X$ contains no nonempty compact
$1$-pseudoconcave subset, and also when $X$ admits a holomorphic map into
some $\CC^N$ of generic rank $n$.
\item[(ii)] Let $v\in\PSH(X)$ satisfy $v\le u$ outside a compact set
but $v>u$ somewhere, and set
\[
 m=\sup_X(v-u),\qquad w=(v-u)^*,\qquad K^*=\{w=m\},
\]
where $^*$ denotes upper semicontinuous regularization. Then $K^*$ is a
nonempty compact $1$-pseudoconcave set which is not MA-positive.
\item[(iii)] In the situation of (ii), $K^*$ contains a minimal compact
$1$-pseudoconcave set $L$. Nonzero positive $\ddc$-closed
currents of bidimension $(1,1)$ supported in $L$ exist, and each of
them has support exactly $L$.
\end{enumerate}
\end{mainthm}

The key analytic step is Lemma~\ref{lem:comparison}, a comparison
principle for locally bounded functions on MA-positive sets, in which the
auxiliary psh function is defined only near the compact set where the
comparison may fail. It may be viewed as a local form of the comparison
framework in \cite{DoHiepComp}. Its proof combines a perturbation inspired
by \cite{DieuUnicity} with a compact-support integration by parts for
mixed Monge--Amp\`ere products. A separate local perturbation shows that
the regularized comparison difference has a pseudoconcave maximum set,
even when $u$ is discontinuous. Shcherbina's criterion, Sibony's duality
and S\l odkowski's minimal-set theorem \cite{Shcherbina,Slodkowski} then
produce the set $L$ and its currents. The two obstructions in (ii) have
different strengths: the blow-up of $\CC^n$ at a point contains a compact
pseudoconcave exceptional divisor, but all its compact subsets are
MA-positive.

Our second main result shows that compact analytic subvarieties are
neither necessary nor sufficient for the failure of the
local-to-global principle.

\begin{mainthm}\label{thmB}
The following assertions hold true.
\begin{enumerate}
\item[(i)] Every simple compact complex torus of dimension $n\ge2$
contains a complete K\"ahler domain $X$ without positive-dimensional
compact analytic subsets on which some bounded smooth psh functions
are locally but not globally maximal.
\item[(ii)] The threefold $(\PP^2\setminus\{x_0\})\times\CC$ is a
complete K\"ahler manifold which contains compact complex curves and on
which every psh function is globally maximal.
\end{enumerate}
\end{mainthm}

The domains in (i) are cut out by a real character of the torus,
which provides a proper pluriharmonic function $h$ whose compact
Levi-flat levels have dense complex leaves. On such manifolds the
maximum principle on compact level bands reduces local and global
maximality to convexity on an interval (Theorem~\ref{thm:dense}), so
that every psh function is locally maximal while the globally maximal
ones are precisely the affine functions of $h$. In complex dimension
two, Weyl's lemma on the leaves and Fourier analysis on the levels
determine all compactly supported positive $\ddc$-closed currents
(Theorems~\ref{thm:support} and~\ref{thm:fourier}). On these surfaces a
compact set is MA-positive exactly when it contains no complete level, so
the two obstructions in Theorem~\ref{thmA}(ii) coincide, in contrast with
the blow-up. The minimal sets $L$ in
Theorem~\ref{thmA} are single levels, whereas comparison defects
require whole bands of levels (Proposition~\ref{prop:band}), so the
converse of Theorem~\ref{thmA} fails. Related Levi-current phenomena
occur on Grauert type surfaces \cite{MST18,BM22,BM24}. In (ii),
psh functions descend through holomorphic submersions whose fibers are
noncompact and carry only constant psh functions (Theorem~\ref{thm:product}). Finally,
Proposition~\ref{prop:singular}, building on Rashkovskii's work
\cite{Rashkovskii1998}, examines how divisorial logarithmic poles
affect local and global maximality.

The paper is organized as follows. Section~\ref{sec:obstruction}
contains the comparison lemma, the basic properties of compactly
supported positive $\ddc$-closed currents and the proof of
Theorem~\ref{thmA}. Section~\ref{sec:models} treats manifolds with
proper pluriharmonic functions, the torus domains and their currents,
the submersions of Theorem~\ref{thmB}(ii), and divisorial poles.

\section{Compact comparison obstructions}\label{sec:obstruction}
We use $1$-pseudoconcavity in the sense of Rothstein
\cite{Shcherbina}. Equivalently, a closed set has the local maximum
property: no smooth strictly psh
function has a strict local maximum on the set
\cite[Section~1]{Slodkowski}. Shcherbina's theorem states
that a compact set admits a smooth strictly psh function on a
neighborhood if and only if it contains no nonempty compact
$1$-pseudoconcave subset. When no such function exists, Sibony's
duality provides a nonzero positive $\ddc$-closed
current of bidimension $(1,1)$ supported on the compact set
\cite[Remark~1.1]{Shcherbina}. We fix a smooth
positive volume form $dV$ on $X$, and we call a compact set
$K\subset X$ \emph{MA-positive} if it has a neighborhood $U$
carrying a bounded psh function $\rho$ such that
\begin{equation}\label{eq:mapositive}
(\ddc\rho)^n\ge f\,dV,\qquad f>0\ \text{almost everywhere on }U.
\end{equation}
Every compact set admitting a smooth strictly psh function on a
neighborhood is MA-positive.

A psh function with positive Monge--Amp\`ere density near the
comparison set allows a compactly supported perturbation and a
Monge--Amp\`ere mass calculation.

\begin{lemma}\label{lem:comparison}
Let $u\in\PSH(X)\cap L^\infty_{\loc}(X)$ be locally maximal, and let $K\subset X$ be an MA-positive compact set. If $v\in\PSH(X)$ satisfies $v\le u$ on $X\setminus K$, then $v\le u$ on $X$.
\end{lemma}
\begin{proof}
Suppose that $v(x_0)>u(x_0)$ for some $x_0\in K$. Let $U$ and $\rho$
be as in \eqref{eq:mapositive}, and choose $K\subset D\cpt U$.
Fix $N>0$ such that $u\ge-N$ on $\overline D$, and put
$\widetilde v=\max\{v,-N\}$ and $M=\sup_{D}\rho$.
For sufficiently small $\varepsilon>0$, the function
$w=\widetilde v+\varepsilon(\rho-M)$ satisfies
$w(x_0)>u(x_0)$ and $w\le u$ on $D\setminus K$.
The function $\widetilde u=\max\{u,w\}$ on $D$, extended by $u$
outside $K$, is locally bounded and psh on $X$.

By local maximality, $(\ddc u)^n=0$. Set
\[
S=\sum_{j=0}^{n-1}(\ddc\widetilde u)^j\wedge(\ddc u)^{n-1-j}.
\]
The multilinear identity for mixed Monge--Amp\`ere products and
integration by parts, since $\widetilde u-u$ has compact support and
$\ddc S=0$, yield
\begin{equation}\label{eq:mass}
\int_X(\ddc\widetilde u)^n
 =\int_X\ddc(\widetilde u-u)\wedge S=0.
\end{equation}
Here the last integral is defined by a cutoff equal to one near $K$.
On the plurifinely open set $G=\{w>u\}$, plurifine locality~\cite{BT87} and
positivity of mixed Monge--Amp\`ere products imply
\begin{equation}\label{eq:positive}
0=\mathbf1_G(\ddc\widetilde u)^n
 =\mathbf1_G(\ddc w)^n
 \ge\varepsilon^n\mathbf1_G(\ddc\rho)^n.
\end{equation}
Since $(\ddc\rho)^n$ has a positive density almost everywhere,
$G$ has Lebesgue measure zero.
If $w(x_0)>u(x_0)$, upper semicontinuity of $u$ makes
$u$ uniformly smaller than $w(x_0)$ on a small coordinate ball
about $x_0$. The submean inequality for $w$ contradicts
the almost-everywhere inequality $w\le u$ on that ball, proving
the comparison principle.
\end{proof}

Compactly supported positive $\ddc$-closed currents annihilate the Levi
forms of smooth psh functions. Their supports are pseudoconcave; these
standard facts about Levi currents \cite{BM22} will be used throughout.

\begin{lemma}\label{lem:currents}
Let $T$ be a nonzero positive $\ddc$-closed current of bidimension
$(1,1)$ with compact support $\Sigma$ on $X$.
For every smooth psh function $w$ defined near $\Sigma$,
\begin{equation}\label{eq:leviann}
 T\wedge\ddc w=0.
\end{equation}
If $h$ is real pluriharmonic near $\Sigma$, then
\begin{equation}\label{eq:directed}
 T\wedge dh\wedge d^ch=0,
\end{equation}
so $T$ is directed by $\ker\partial h$ wherever $dh\ne0$.
Moreover, $\Sigma$ is a compact $1$-pseudoconcave set.
\end{lemma}
\begin{proof}
Choose a smooth cutoff $\theta$ supported in the domain of $w$ and
equal to one near $\Sigma$. The positive measure $T\wedge\ddc w$ has total mass
\[
 \langle T,\theta\ddc w\rangle
 =\langle\ddc T,\theta w\rangle=0.
\]
Applying this to $w=h^2$ proves \eqref{eq:directed}.
In coordinates with $h=\textup{Re}\,z_1$, positivity and the
Cauchy--Schwarz inequality then eliminate the normal and mixed
coefficients of $T$, proving directedness.

If the local maximum property failed on $\Sigma$, there would be
$x_0\in\Sigma$, a closed coordinate ball $\overline B$ about $x_0$,
and a smooth strictly psh function $\rho$ near $\overline B$ such that
$\rho(x_0)=0$ and $\rho<0$ on $\Sigma\cap\partial B$.
Choose $\eta>0$ and a neighborhood $W$ of $\Sigma\cap\partial B$
on which $\rho<-3\eta$. Let $m_\eta$ be a smooth regularized maximum
which equals $\max\{s,t\}$ when $|s-t|\ge\eta$ and is convex and
nondecreasing in both variables. The function
\[
 w=m_\eta(\rho,-2\eta)+2\eta\quad\text{on }B,
 \qquad w=0\quad\text{on }W\cup(X\setminus\overline B)
\]
is smooth and psh on a neighborhood of $\Sigma$, and equals
$\rho+2\eta$ near $x_0$.
By \eqref{eq:leviann}, $T\wedge\ddc w=0$, whereas near $x_0$
this measure dominates a positive multiple of the trace measure of
$T$. Since $x_0\in\supp T$, this is a contradiction.
\end{proof}

The upper semicontinuous regularization of the comparison difference
locates a single obstruction shared by all its positive superlevels.

\begin{theorem}\label{thm:obstruction}
Let $u\in\PSH(X)\cap L^\infty_{\loc}(X)$ be locally maximal, and let
$v\in\PSH(X)$ satisfy $v\le u$ outside a compact set but $v>u$
somewhere. Set
\begin{equation}\label{eq:maxset}
 m=\sup_X(v-u),\qquad w=(v-u)^*,\qquad K^*=\{w=m\}.
\end{equation}
Then the following assertions hold true.
\begin{enumerate}
\item[(i)] $K^*$ is a nonempty compact $1$-pseudoconcave set which is
not MA-positive.
\item[(ii)] $K^*$ contains a minimal compact $1$-pseudoconcave set $L$. Nonzero positive $\ddc$-closed currents of bidimension $(1,1)$
supported in $L$ exist, each of them has support exactly $L$, and one of
them spans an extremal ray in the cone of positive $\ddc$-closed currents
of bidimension $(1,1)$ on $X$.
\item[(iii)] For every $0<\delta<m$, the compact set
\begin{equation}\label{eq:superlevel}
 F_\delta=\overline{\{v>u+\delta\}}
\end{equation}
contains $L$ and is not MA-positive.
\item[(iv)] If $u$ is continuous, then $K^*=\{v-u=m\}$.
\end{enumerate}
\end{theorem}
\begin{proof}
We first prove (i) and (iii). The number $m$ is finite and positive because $u$ is locally bounded
and $v\le u$ off a compact set $K$. The sets $F_\delta$ are nonempty
compact subsets of $K$, decrease as $\delta\uparrow m$, and satisfy
\[
 K^*=\bigcap_{0<\delta<m}F_\delta.
\]
Thus $K^*$ is nonempty and compact. The inequality $v-\delta\le u$
holds off $F_\delta$ and fails somewhere, so Lemma~\ref{lem:comparison}
shows that $F_\delta$ is not MA-positive.
If $K^*$ were MA-positive on a neighborhood $U$, compactness and the
decreasing-intersection identity would imply $F_\delta\subset U$ for
$\delta$ sufficiently close to $m$. The same auxiliary function would
make $F_\delta$ MA-positive, a contradiction.

We next verify the local maximum property of $K^*$.
If it failed, choose $x_0\in K^*$, a small coordinate ball $B$ about
$x_0$ inside a neighborhood on which $u$ is maximal, and a smooth
strictly psh function $\rho$ near $\overline B$ such that
$\rho(x_0)=0$ and $\rho<0$ on $K^*\cap\partial B$.
Since $w$ is upper semicontinuous and $w\le m$, for sufficiently
small $\varepsilon>0$ we have
\[
 \max_{\partial B}(w+\varepsilon\rho)<m.
\]
Choose $\eta>0$ such that
$\widetilde v=v+\varepsilon\rho-m+\eta<u$ in a collar of $\partial B$.
This follows from $v-u\le w$ and the strict boundary inequality.
The definition $w(x_0)=m$ provides points arbitrarily close to $x_0$
at which $v-u+\varepsilon\rho>m-\eta$, so $\widetilde v>u$ somewhere
in $B$. Hence $\max\{u,\widetilde v\}$ on $B$, extended by $u$, contradicts
the local maximality of $u$. This proves pseudoconcavity.

For (ii), by \cite[Proposition~1.5]{Slodkowski}, $K^*$ contains
a minimal compact $1$-pseudoconcave set $L$. Sibony's duality
\cite[Remark~1.1]{Shcherbina} supplies a nonzero positive
$\ddc$-closed current supported in $L$.
For every such current $T$, Lemma~\ref{lem:currents} and minimality
imply $\supp T=L$.
Normalize these currents to have mass one with respect to a Hermitian
metric near $L$. They form a nonempty weakly compact convex set, and an
extreme point spans an extremal ray of the cone in (ii), since every
positive summand of it has support contained in $L$. Finally, $L\subset K^*\subset F_\delta$
for every $\delta$, and continuity of $u$ makes $v-u$ upper
semicontinuous, proving the last assertion.
\end{proof}

By \cite[Corollary~1.7]{Slodkowski}, every psh function defined near
the minimal set $L$ is constant on $L$. Thus $L$ behaves like a compact
complex curve, which explains why it obstructs global comparison.

The obstruction theorem yields a local-to-global principle when
compact pseudoconcave sets are absent.

\begin{corollary}\label{cor:stein}
Assume that every compact subset of $X$ is MA-positive. Then every
locally bounded locally maximal psh function on $X$ is globally
maximal. The assumption holds in each of the following cases:
\begin{enumerate}
\item[(i)] $X$ contains no nonempty compact $1$-pseudoconcave subset;
\item[(ii)] $X$ admits a holomorphic map $F=(F_1,\dots,F_N):X\to\CC^N$
of generic rank $n$.
\end{enumerate}
\end{corollary}
\begin{proof}
The first assertion is Lemma~\ref{lem:comparison}. In case (i),
Shcherbina's theorem provides a smooth strictly psh function near
every compact subset of $X$. In case (ii), the smooth psh function
$\rho=\sum_j|F_j|^2$ is bounded near every compact set, and
\[
\ddc\rho=\frac{\ii}{\pi}\sum_{j=1}^N\partial F_j\wedge\overline{\partial F_j}
\]
has rank equal to the rank of $dF$ at each point. Hence
$(\ddc\rho)^n>0$ outside the analytic set where $\operatorname{rank}dF<n$,
which is a proper analytic subset of the connected manifold $X$ and
therefore has Lebesgue measure zero.
\end{proof}

Theorem~\ref{thmA} follows from Corollary~\ref{cor:stein} and
Theorem~\ref{thm:obstruction}.

Stein manifolds satisfy both (i) and (ii), and for them the corollary
is well known. The blow-up $\pi:X\to\CC^n$ of the origin satisfies
(ii) but not (i): its exceptional divisor is a compact complex
submanifold, hence $1$-pseudoconcave by the maximum principle along
it. Thus MA-positivity is strictly weaker than the existence of a
smooth strictly psh function nearby, and failure of MA-positivity is a
stronger obstruction than the presence of a compact pseudoconcave subset.
In particular, compact pseudoconcave sets do not by themselves force
failure of global maximality; Proposition~\ref{prop:band} and
Theorem~\ref{thm:product} below provide further reasons for this.
We do not know a dual characterization of compact sets which are not
MA-positive, analogous to the description of compact sets without
smooth strictly psh neighborhoods through Sibony's duality.

\section{Complete K\"ahler models and directed currents}\label{sec:models}
Convex functions of a pluriharmonic function have Levi form of
rank at most one, which accounts for local maximality in dimension
at least two. Properness allows a chord replacement on an interval
to produce a competitor that differs only on a compact set.

\begin{proposition}\label{prop:convex}
Let $X$ be a connected noncompact complex manifold of dimension $n\ge2$, and let $h:X\to I$ be a proper surjective real pluriharmonic function onto an open interval $I$. For every finite convex function $g:I\to\RR$, the composition $g\circ h$ is continuous, psh and locally maximal. Moreover,
\begin{equation}\label{eq:affine}
g\circ h\in\MPSH(X)\quad\text{if and only if}\quad g\text{ is affine}.
\end{equation}
\end{proposition}
\begin{proof}
A finite convex function on $I$ is continuous and admits smooth convex approximations converging uniformly on compact subintervals. For smooth convex $g$, the identity $\ddc h=0$ implies
\begin{equation}\label{eq:ddccomposition}
\ddc(g\circ h)=g''(h)\,dh\wedge d^ch,
\end{equation}
a semipositive form of rank at most one. Its $n$th power vanishes since $n\ge2$. Continuity of the Monge--Amp\`ere operator under locally uniform
convergence implies $(\ddc(g\circ h))^n=0$ for arbitrary finite convex $g$, hence local maximality.

If $g$ is affine, then $g\circ h$ is pluriharmonic. Suppose $v\in\PSH(X)$ satisfies $v\le g\circ h$ off a compact set $K$. Choose a relatively compact domain $D\supset K$. The psh function $v-g\circ h$ is nonpositive on $\partial D\subset X\setminus K$, hence on $D$ by the maximum principle. Thus $g\circ h$ is globally maximal.

If $g$ is not affine, choose $c<d$ in $I$ such that the chord $\ell$ through $(c,g(c))$ and $(d,g(d))$ lies strictly above $g$ at some point of $(c,d)$. The function agreeing with $\ell$ on $[c,d]$ and with $g$ outside $[c,d]$ is convex by monotonicity of secant slopes. Its composition with $h$ is psh, strictly exceeds $g\circ h$ somewhere, and agrees with it outside the compact set $h^{-1}([c,d])$,
so $g\circ h$ is not globally maximal.
\end{proof}

At a regular value $t$, the level $M_t=h^{-1}(t)$ is Levi-flat.
Its complex leaves are obtained by continuing the local plaques
$\{F=\mathrm{constant}\}$, where $h=\textup{Re}\,F$ for a local
holomorphic function $F$ with $dF\ne0$.
Factorization through a pluriharmonic function is familiar on Grauert
type surfaces \cite{MST18}. The following level-maximum argument works
in every dimension and only requires dense leaves at a dense set of
regular values.

\begin{theorem}\label{thm:dense}
Let $X$ be a connected complex manifold of dimension $n\ge2$, and let $h:X\to I$ be a proper surjective real pluriharmonic function onto an open interval $I$. Assume that there is a dense subset $J\subset I$ consisting of regular
values such that every complex leaf of $M_t=h^{-1}(t)$ is dense
in $M_t$ for each $t\in J$.
Then
\begin{equation}\label{eq:pshclass}
\begin{gathered}
\PSH(X)=\{g\circ h:g\text{ finite and convex on }I\},\\
\MPSH(X)=\{c_1h+c_2:c_1,c_2\in\RR\}.
\end{gathered}
\end{equation}
Thus every psh function on $X$ is continuous and locally maximal; global maximality is equivalent to pluriharmonicity.
\end{theorem}
\begin{proof}
Fix $w\in\PSH(X)$ and define
\[
 g(t)=\max_{h^{-1}(t)}w.
\]
We first show that $g$ is finite. Every level of $h$ contains a
regular point. Indeed, the critical locus is a proper analytic set
of real codimension at least two. Its complement in a sufficiently
small coordinate ball is path-connected, and $h-t$ assumes both
signs near every point of its zero set. A path joining points of
opposite signs outside the critical locus therefore meets the level
at a regular point. Near such a point, choose coordinates with
$h=t+\textup{Re}\,z_1$. If $w$ were identically $-\infty$ on
$M_t$, each subharmonic restriction in the $z_1$ direction would be
$-\infty$ on a nonpolar real segment, hence identically $-\infty$.
This would force $w\equiv-\infty$ on $X$. Thus $g(t)>-\infty$;
compactness of the level implies $g(t)<+\infty$.

For $c<d$ in $I$, let $\ell$ be the affine function with
$\ell(c)=g(c)$ and $\ell(d)=g(d)$. The psh function
$w-\ell\circ h$ is nonpositive on the boundary of the relatively
compact open set $h^{-1}((c,d))$. The maximum principle yields
\[
 g(t)\le\frac{d-t}{d-c}g(c)+\frac{t-c}{d-c}g(d),
 \qquad c<t<d.
\]
Consequently $g$ is finite, convex and continuous on $I$.

For $t\in J$, choose a maximum point of $w$ on $M_t$.
The restriction of $w$ to its complex leaf attains the finite
maximum $g(t)$ and is therefore constant. Density of that leaf
and upper semicontinuity imply $w=g(t)$ on $M_t$.
The nonconstant pluriharmonic function $h$ is an open map.
Hence every $x\in X$ is a limit of points $x_j\in h^{-1}(J)$.
Upper semicontinuity and continuity of $g$ imply
\[
 w(x)\ge\limsup_{j\to\infty}w(x_j)
       =g(h(x)).
\]
The reverse inequality follows from the definition of $g$.
Thus $w=g\circ h$ everywhere. Proposition~\ref{prop:convex}
proves the converse inclusion and the maximality assertion.
Finally, if $g\circ h$ is pluriharmonic, applying the same
classification to its negative shows that both $g$ and $-g$ are
convex, so $g$ is affine.
\end{proof}

\begin{remark}\label{rem:density}
The density hypothesis cannot be dropped. On $X=E\times\{r<|z|<R\}$,
with $E$ an elliptic curve, the levels of $h(x,z)=\log|z|$ have
compact leaves $E\times\{z\}$, and the pluriharmonic function
$\textup{Re}\,z$ does not factor through $h$.
\end{remark}

On surfaces, the support of a directed positive $\ddc$-closed current
is a union of leaves. The resulting full-level support property is
standard for Grauert type surfaces; see \cite[Theorem~6.1]{BM22}.
We give a short local proof for the pluriharmonic submersions used here.

\begin{theorem}\label{thm:support}
Let $h:X\to I$ be a proper surjective real pluriharmonic
submersion from a connected complex surface onto an open interval. Suppose that every complex leaf
is dense in its level. Then every nonzero positive $\ddc$-closed current $T$ of bidimension $(1,1)$ with compact support has
\begin{equation}\label{eq:support}
\supp T=h^{-1}(S)
\end{equation}
for a nonempty compact subset $S\subset I$. Conversely, every level $M_t$ is the support of a positive closed
$(1,1)$-current. Thus, for a compact set $K\subset X$, the following assertions are equivalent:
\begin{enumerate}
\item[(i)] $K$ contains a nonempty compact $1$-pseudoconcave subset;
\item[(ii)] $K$ supports a nonzero positive $\ddc$-closed current of bidimension $(1,1)$;
\item[(iii)] $K$ contains a whole level $M_t$ for some $t\in I$;
\item[(iv)] $K$ is not MA-positive.
\end{enumerate}
\end{theorem}
\begin{proof}
By Lemma~\ref{lem:currents}, $T$ is directed by $\ker\partial h$.
In a chart $U\simeq\mathbb D_1\times\mathbb D_2$ with
$h=\textup{Re}\,z_1+\mathrm{constant}$, positivity therefore implies
\[
 T=\ii\mu\,dz_1\wedge d\bar z_1,
 \qquad \frac{\partial^2\mu}{\partial z_2\,\partial\bar z_2}=0
\]
for a positive measure $\mu$. For each nonnegative
$w\in C_c^\infty(\mathbb D_1)$, the measure
\[
 \langle\mu_w,\widetilde w\rangle
 =\int_Uw(z_1)\widetilde w(z_2)\,d\mu
\]
has a nonnegative harmonic density by Weyl's lemma. This density is
either zero or everywhere positive. Applying this to test functions
supported near any given $z_1$ shows that a plaque meeting
$\supp T$ is contained in $\supp T$.
Thus the support is a union of complete leaves and, by density and
closedness, of whole levels. Equation~\eqref{eq:support} follows with
$S=h(\supp T)$.

Conversely, the positive closed current
$Q_t=\ddc\max\{h-t,0\}$ has support $M_t$: locally it is a positive
multiple of $\delta_0(\textup{Re}\,z_1)\,\ii dz_1\wedge d\bar z_1$.
For the equivalence, (i) implies (ii) by Shcherbina's theorem and
Sibony's duality, while the support description proves (ii) implies
(iii). If $K\supset M_t$ were MA-positive on a neighborhood $U$, then
properness would yield a band $h^{-1}((t-\varepsilon,t+\varepsilon))\subset U$.
Each connected component $C$ of this band is mapped by $h$ properly
and openly, hence onto $(t-\varepsilon,t+\varepsilon)$, and every
level of $h|_C$ is a full level of $h$ because the leaves are dense.
Theorem~\ref{thm:dense}, applied on $C$, writes the auxiliary function
$\rho$ as a convex function of $h$, so $(\ddc\rho)^2=0$ on $C$, a
contradiction. Hence (iii) implies (iv).
Finally, a set which is not MA-positive admits no smooth strictly psh
function nearby, and Shcherbina's theorem shows that (iv) implies (i).
\end{proof}

Comparison defects require bands, whereas by Theorem~\ref{thm:support} positive currents, and the failure of MA-positivity, may occur on single levels.

\begin{proposition}\label{prop:band}
Under the hypotheses of Theorem~\ref{thm:dense}, let $K\cpt X$ and
$u=g\circ h$, with $g:I\to\RR$ finite and convex.
There exists $v\in\PSH(X)$ such that $v\le u$ on $X\setminus K$
but $v>u$ somewhere if and only if, for some $c<d$ in $I$,
\[
 h^{-1}([c,d])\subset K,
 \qquad g|_{[c,d]}\text{ is not affine}.
\]
In particular, $K$ supports a comparison defect for some locally
maximal psh function if and only if it contains a complete band
$h^{-1}([c,d])$, where $c<d$.
\end{proposition}
\begin{proof}
By Theorem~\ref{thm:dense}, every competitor has the form
$v=\widetilde g\circ h$, where $\widetilde g$ is finite and convex. If $\widetilde g>g$
somewhere, let $(c,d)$ be a component of $\{\widetilde g>g\}$.
Since every corresponding level lies in $K$ and $h(K)$ is compact,
we have $c,d\in I$. Every point of $h^{-1}(\{c,d\})$
is approached by levels in $(c,d)$, since a nonconstant
pluriharmonic function has no local extremum. Hence
$h^{-1}([c,d])\subset K$. At the endpoints $\widetilde g=g$; convexity
of $\widetilde g$ shows that $g$ cannot be affine on $[c,d]$.

Conversely, let $\ell_{c,d}$ be the chord joining $(c,g(c))$ and
$(d,g(d))$. The function equal to $\ell_{c,d}$ on $[c,d]$ and to
$g$ elsewhere is convex, by the secant-slope inequalities.
Its composition with $h$ is an admissible competitor and exceeds
$u$ somewhere. Taking $g(t)=t^2$ proves the last assertion.
\end{proof}

A primitive lattice character on a simple complex torus produces
compact Levi-flat levels with dense leaves.

\begin{theorem}\label{thm:abelian}
Let $A$ be any simple compact complex torus of dimension $n\ge2$.
There exists a connected smoothly bounded domain $X\subset A$ such that:
\begin{enumerate}
\item[(i)] $X$ is a complete K\"ahler manifold with no positive-dimensional compact analytic subset;
\item[(ii)] $X$ admits a proper real pluriharmonic submersion $h:X\to(a,b)$ whose levels are connected compact real-analytic Levi-flat hypersurfaces with dense complex leaves;
\item[(iii)] every psh function on $X$ is locally maximal, and
\begin{equation}\label{eq:abelianclass}
\begin{gathered}
\PSH(X)=\{g\circ h:g\text{ finite and convex on }(a,b)\},\\
\MPSH(X)=\{c_1h+c_2:c_1,c_2\in\RR\}.
\end{gathered}
\end{equation}
\item[(iv)] $X$ admits bounded smooth psh functions, and smooth psh exhaustions, that are locally but not globally maximal.
\end{enumerate}
Every level $M_t=h^{-1}(t)$ supports a nonzero positive \emph{closed} current of bidimension $(1,1)$ having full support on $M_t$.
\end{theorem}
\begin{proof}
Let $A=\CC^n/\Lambda$. Fix a surjective group homomorphism $\lambda:\Lambda\to\ZZ$. Its real linear extension to $\CC^n$ is of the form $\textup{Re}L$ for a complex linear functional $L$, and hence
\[
f:A\longrightarrow\RR/\ZZ,\qquad
f([z])=\textup{Re}L(z)\pmod{\ZZ}
\]
is a surjective homomorphism of real tori. Because $\lambda$ is primitive, $H=\ker f$ is connected. Choose $0<a<b<1$ and put
\[
X=f^{-1}((a,b)\pmod{\ZZ}).
\]
The map $f|_X$ has a proper pluriharmonic lift $h:X\to(a,b)$, and each $M_t$ is a translate of the connected compact real torus $H$ of dimension $2n-1$. The boundary of $X$ is the disjoint union of the two smooth levels $M_a$ and $M_b$ in $A$.

Put $V=\ker L$. The leaf through each point of $M_t$ is the image of a translate of $V\simeq\CC^{n-1}$. The closure of this image is a connected real subtorus of $H$ of real dimension at least $2n-2$. If it were a proper subtorus of $H$, its dimension would be $2n-2$ and its Lie algebra would be $V$. It would therefore be a positive-dimensional proper complex subtorus of $A$, contradicting simplicity. Hence the complex leaves are dense, and Theorem~\ref{thm:dense} yields (iii).

Suppose an irreducible positive-dimensional compact analytic set $Y$
lies in $X$. The maximum principle on a resolution of $Y$ makes $h$
constant on $Y$. Thus $\Sigma=Y-Y$ is a positive-dimensional compact
irreducible analytic subset of $H$ containing $0$.
For $m\ge1$, let $\Sigma_m=\Sigma+\dots+\Sigma$ be the set of sums of
$m$ elements of $\Sigma$. These sets are compact, irreducible, analytic and
increasing, by the proper mapping theorem and since $0\in\Sigma$. Stabilization of their dimensions implies
$\Sigma_m=\Sigma_{m+1}$ for some $m$.
For $z\in\Sigma$, the translate $\Sigma_m+z\subset\Sigma_m$ has the
same dimension, hence equals $\Sigma_m$.
Iterating, we obtain $\Sigma_m+y=\Sigma_m$ for every $y\in\Sigma_m$.
Since $0\in\Sigma_m$, these identities show that $\Sigma_m$ is closed
under addition and, as $0\in\Sigma_m+y$, under inversion. Hence
$\Sigma_m$ is a connected compact complex subgroup of $A$, that is, a
complex subtorus. It is positive-dimensional
and contained in $H\ne A$, contradicting simplicity.

To construct a complete metric, let $\omega_0$ be a flat K\"ahler form on $A$ and set
\begin{equation}\label{eq:complete}
G(t)=-\log(t-a)-\log(b-t),\qquad
\omega=\omega_0+\ddc(G\circ h)
=\omega_0+G''(h)\,dh\wedge d^ch.
\end{equation}
The form $\omega$ is K\"ahler, and
\[
G''(t)=\frac{1}{(t-a)^2}+\frac{1}{(b-t)^2}.
\]
For a real tangent vector $Y$ we have $|dh(Y)|\le2|\partial h(Y)|$,
so the $\omega$-length of a curve $\gamma$ is at least a positive
multiple of $\int\sqrt{G''(h\circ\gamma)}\,|(h\circ\gamma)'|$.
Since $\sqrt{G''(t)}$ dominates $\frac{1}{t-a}$ and $\frac{1}{b-t}$,
which are not integrable at the endpoints, every bounded subset of
$(X,\omega)$ lies in $h^{-1}([c,d])$ for some $a<c<d<b$. This set is
compact because $h$ is proper, and the Hopf--Rinow theorem shows that
$\omega$ is complete.

For $t_0=\frac{a+b}{2}$, the bounded smooth psh function
$u=(h-t_0)^2$ is locally but not globally maximal by
Proposition~\ref{prop:convex}. The same holds for the smooth psh
exhaustion $g\circ h$, where $g(t)=\frac{1}{t-a}+\frac{1}{b-t}$.
Proposition~\ref{prop:band} describes exactly which compact sets
can support these comparison defects.

Finally, choose $0\ne\xi\in V$ and regard it as a translation-invariant holomorphic vector field on $A$. Denote by $\sigma_t$ normalized Haar measure on $M_t$ and define, for test $(1,1)$-forms $\alpha$,
\begin{equation}\label{eq:haar}
\langle T_t,\alpha\rangle
=\int_{M_t}\frac1{\ii}\alpha(\xi,\bar\xi)\,d\sigma_t.
\end{equation}
The current is positive and has full support on $M_t$.
The real and imaginary flows generated by $\xi$ preserve Haar
measure. For a compactly supported smooth complex $1$-form $\beta$,
\[
 d\beta(\xi,\bar\xi)
 =\xi\bigl(\beta(\bar\xi)\bigr)
  -\bar\xi\bigl(\beta(\xi)\bigr),
\]
since the invariant vector fields commute. Both terms integrate
to zero against $\sigma_t$, so $T_t$ is $d$-closed.
\end{proof}

Simple projective abelian varieties exist in every dimension
\cite{BirkenhakeLange}, so Theorem~\ref{thm:abelian} is not vacuous.

For abelian surfaces, Fourier analysis of the leafwise equation
classifies compactly supported positive $\ddc$-closed currents.

\begin{theorem}\label{thm:fourier}
Let $X$ be the two-dimensional torus domain in
Theorem~\ref{thm:abelian}, and let $T_t$ be the Haar currents
from \eqref{eq:haar}. Every compactly supported positive
$\ddc$-closed current $T$ of bidimension $(1,1)$ on $X$
has a unique representation
\begin{equation}\label{eq:disintegrate}
 T=\int_{(a,b)}T_t\,d\nu(t),
\end{equation}
where $\nu$ is a finite positive Borel measure with compact
support in $(a,b)$. Conversely, every such measure defines
one of these currents. Moreover,
\begin{equation}\label{eq:disintsupport}
 dT=0,\qquad \supp T=h^{-1}(\supp\nu).
\end{equation}
\end{theorem}
\begin{proof}
By Lemma~\ref{lem:currents}, every current in the theorem
is directed by $\ker\partial h$. We first consider any positive
$\ddc$-closed current directed by this line, without a
compact-support assumption. Since the line is spanned by $\xi$,
the current has the form
\begin{equation}\label{eq:mu}
 \langle T,\alpha\rangle
 =\int_X\frac1{\ii}\alpha(\xi,\bar\xi)\,d\mu
\end{equation}
for a unique positive locally finite measure $\mu$ on $X$.
The equation $\ddc T=0$ is equivalent to
\begin{equation}\label{eq:xi}
 \int_X\xi\bar\xi(w)\,d\mu=0
 \qquad(w\in C_c^\infty(X)).
\end{equation}

Choose $\gamma\in\Lambda$ with $\lambda(\gamma)=1$. Translation in the real direction $\RR\gamma$ splits the real tori $A\simeq H\times(\RR/\ZZ)$, $H=\ker f$, and yields a real-analytic trivialization $X\simeq(a,b)\times H$. The real and imaginary parts $P,Q$ of $\xi$ are constant vector fields on the three-torus $H$, and $\xi\bar\xi$ is a positive scalar multiple of $P^2+Q^2$ on functions. For a character $e_k$ of $H$, $k\in\ZZ^3$, we have
\[
Pe_k=2\pi\ii\langle k,P\rangle e_k,
\qquad Qe_k=2\pi\ii\langle k,Q\rangle e_k.
\]
Then $(P^2+Q^2)e_k=-c_ke_k$, where
\[
c_k=4\pi^2\bigl(\langle k,P\rangle^2+\langle k,Q\rangle^2\bigr).
\]
If $k\ne0$ then $c_k>0$; otherwise $e_k$ would be invariant under both flows, whose jointly generated subgroup is dense in $H$ by Theorem~\ref{thm:abelian}. Testing \eqref{eq:xi} with $w(t,y)=g_1(t)e_k(y)$, where $g_1\in C_c^\infty((a,b))$, we obtain
\[
\int_X g_1(t)e_k(y)\,d\mu(t,y)=0\qquad(k\ne0).
\]
All nonconstant fiberwise Fourier coefficients of $\mu$ therefore vanish. As trigonometric polynomials are dense in $C(H)$,
\[
d\mu(t,y)=d\nu(t)\,d\sigma_t(y),
\qquad\nu=h_*\mu.
\]
Substitution in \eqref{eq:mu} proves \eqref{eq:disintegrate}
and its uniqueness. Haar measure has full support on each level, so
\eqref{eq:disintsupport} follows. Conversely, integrating the
closed currents $T_t$ against a locally finite $\nu$ defines a
positive closed current: each test form meets only a compact
interval of levels. When $T$ has compact support, the support formula and properness
of $h$ imply that $\nu$ has compact support.
\end{proof}

In particular, the extremal rays of the cone of compactly supported
positive $\ddc$-closed currents on this surface are generated by the
Haar currents $T_t$. Since the minimal sets $L$ in
Theorem~\ref{thm:obstruction} are single levels by
Theorem~\ref{thm:support}, the currents they carry are multiples of
Haar currents.

\begin{remark}\label{rem:higherdim}
The full-level support conclusion is specific to surfaces. For $n\ge3$,
choose a second integer character $\lambda_2$ such that
$(\lambda,\lambda_2)$ is a primitive pair with linearly independent
complex-linear extensions $L,L_2$; this is possible as soon as $n\ge2$,
while $n\ge3$ guarantees some $0\ne\zeta\in\ker L\cap\ker L_2$. Replacing
$\xi$ and $M_t$ in \eqref{eq:haar} by $\zeta$ and a joint level
$H_{t,s}\subsetneq M_t$ of the two characters, we obtain a positive
closed current with support exactly $H_{t,s}$.
\end{remark}

To illustrate Theorem~\ref{thm:obstruction}, let $u=(h-t_0)^2$ on the two-dimensional torus domain, where $t_0=\frac{a+b}{2}$, and $v=\max\{u,r^2\}$ with $0<r<\min\{t_0-a,b-t_0\}$. For $0<\delta<r^2$ we have $F_\delta=\{|h-t_0|\le\sqrt{r^2-\delta}\}$, the maximum set $K^*$ is the level $M_{t_0}$, which is also the minimal set $L$, and the Haar current $T_{t_0}$ realizes assertion~(ii) of Theorem~\ref{thm:obstruction}. For $n\ge3$, Remark~\ref{rem:higherdim} and Lemma~\ref{lem:currents} show that $L$ is strictly smaller than $M_{t_0}$.

Holomorphic submersions whose fibers are noncompact and carry only
constant psh functions satisfy a direct comparison property, since every
fiber meets the complement of every compact set.

\begin{theorem}\label{thm:product}
Let $\pi:X\to S$ be a surjective holomorphic submersion between connected
complex manifolds. Assume that every fiber is connected and noncompact,
and that every psh function on each fiber is constant. Then
\begin{equation}\label{eq:productpsh}
 \PSH(X)=\{g\circ\pi:g\in\PSH(S)\}=\MPSH(X).
\end{equation}
In particular, the threefold
$X=(\PP^2\setminus\{x_0\})\times\CC$ admits a complete K\"ahler metric,
contains compact complex curves, and has nonconstant smooth psh functions,
while every psh function on $X$ is globally maximal.
\end{theorem}
\begin{proof}
The restriction of $w\in\PSH(X)$ to each fiber is constant or
identically $-\infty$, so $w=g\circ\pi$ for a function $g$ on $S$.
Local holomorphic sections of the submersion show that $g$ is psh,
since $g\equiv-\infty$ on an open set would force $w\equiv-\infty$
on its preimage.
Conversely, every psh pullback is psh.
Suppose $u=g\circ\pi$ and $v=\widetilde g\circ\pi$ satisfy $v\le u$
outside a compact set $K$. Every fiber meets $X\setminus K$, so
$\widetilde g(s)\le g(s)$ for every $s\in S$, proving global maximality.

For $Y=\PP^2\setminus\{x_0\}$, every point lies on a projective line
avoiding $x_0$. A psh function is constant, possibly $-\infty$, on each
such line, and any two of these lines meet in $Y$. Thus every psh function on $Y$ is constant,
and the projection $Y\times\CC\to\CC$ satisfies the hypotheses.
The function $|z|^2$ is a nonconstant smooth psh function on this product,
which contains a compact curve $L\times\{0\}$ for every projective line
$L$ avoiding $x_0$. Demailly's theorem \cite[Theorem~1.5]{Demailly82}
provides a complete K\"ahler metric on $Y$, since it is the complement
of an analytic set in the compact K\"ahler manifold $\PP^2$.
Its product with the flat metric on $\CC$ is complete.
\end{proof}

Theorem~\ref{thmB} follows from Theorem~\ref{thm:abelian}, assertions (i) and (iv), and from Theorem~\ref{thm:product}.

For $S=\mathbb D^2$, let
$\pi_S:(\PP^2\setminus\{x_0\})\times S\to S$ be the projection and set
$g(z)=\log\sqrt{|z_1|^2+|z_2|^2}$.
Theorem~\ref{thm:product} makes $g\circ\pi_S$ globally maximal, with a
pole set of codimension two, whereas $(\ddc g)^2$ is the Dirac mass at
the origin and $g$ is not maximal on $S$.

\medskip
The next result extends Rashkovskii's local addition theorem
\cite[Theorem~1]{Rashkovskii1998} to global maximality with
noncompact divisorial poles.

\begin{proposition}\label{prop:singular}
Let $X$ be a connected noncompact complex manifold of dimension
$n\ge2$, and let $u_D\in\PSH(X)$ satisfy
\[
 \ddc u_D=[D],\qquad D=\sum_j c_jH_j,
\]
where $D$ is a locally finite effective real divisor. Then the
following assertions hold true for every $u\in\PSH(X)$.
\begin{enumerate}
\item[(i)] The function $u$ is locally maximal if and only if
$u+u_D$ is locally maximal.
\item[(ii)] If $u+u_D\in\MPSH(X)$, then $u\in\MPSH(X)$.
\item[(iii)] If each irreducible component $H_j$ of $D$ is
noncompact and $u\in\MPSH(X)$, then $u+u_D\in\MPSH(X)$.
\end{enumerate}
\end{proposition}
\begin{proof}
We first prove (iii). Let $v\in\PSH(X)$ satisfy $v\le u+u_D$
outside a compact set $K$.
For each $H_j$, choose a smooth point of $H_j\setminus K$
which lies on no other component. Near that point,
$u_D=c_j\log|f_j|+h_j$, where $f_j$ defines $H_j$ and $h_j$
is pluriharmonic. Monotonicity of Lelong numbers implies
$\nu(v,x)\ge\nu(u_D,x)=c_j$ for all $x$ in a nonempty open subset
of $H_j$. By Siu's analyticity theorem \cite{Siu1974}, the set
$\{x:\nu(v,x)\ge c_j\}$ is analytic; since it contains a nonempty
open subset of the irreducible set $H_j$, it contains $H_j$. Hence
the generic Lelong number of $v$ along $H_j$ is at least $c_j$, and
Siu's divisorial decomposition \cite{Siu1974} shows that
\[
 \ddc v\ge\sum_jc_j[H_j]=\ddc u_D.
\]
Thus $v-u_D$ has a psh representative $w$ on $X$.
The psh functions $v$ and $w+u_D$ agree almost everywhere,
hence everywhere. On $X\setminus(K\cup\supp D)$,
$w\le u$. Since inequalities between psh functions holding almost
everywhere hold everywhere, $w\le u$ on $X\setminus K$.
Global maximality of $u$ yields $w\le u$ on $X$, so
$v\le u+u_D$ on $X$.

For (ii), suppose that $u+u_D$ is globally maximal.
If $v\le u$ outside $K$, applying global maximality to
$v+u_D$ shows that $v+u_D\le u+u_D$ everywhere.
Cancellation off $\supp D$, followed by the same almost-everywhere
comparison, shows that $v\le u$ everywhere.

To prove (i), work on sufficiently small coordinate balls.
Every irreducible hypersurface in such a ball is noncompact,
as a Stein domain contains no positive-dimensional compact
analytic subsets. The restriction of a maximal psh function to
a smaller domain remains maximal. Applying (ii) and (iii) on these
balls shows that $u$ and $u+u_D$ have the same local maximality
property.
\end{proof}

Assertion~(ii) requires no assumption on $D$; for a related local
criterion in dimension two, see \cite[Theorem~1.5]{DoHiepComp}.

In particular, let $F=g\widetilde F$, where
$0\ne g\in\mathcal O(X)$ and $\widetilde F:X\to\CC^m$
is a nonzero holomorphic mapping. If every irreducible component
of $\{g=0\}$ is noncompact, then $\log\|F\|$ is globally
maximal exactly when $\log\|\widetilde F\|$ is. This is the global counterpart of
the common-factor phenomenon in \cite{Rashkovskii1998}.
For $g(z)=z_1$ on $\CC^n$ with $n\ge2$, the function
$\log|g|$ is globally maximal although it is not locally bounded.

The noncompactness assumption in (iii) cannot simply be dropped.
On $Z\times\mathbb D$, where $Z$ is compact, connected and
positive-dimensional, the function
$u_D(x,z)=\log|z|$ satisfies $\ddc u_D=[Z\times\{0\}]$.
The constant function $u=0$ is globally maximal, whereas
$u_D$ is not: for $0<r<1$, the function
$\max\{u_D,\log r\}$ exceeds $u_D$ and agrees with
it outside the compact set $Z\times\{|z|\le r\}$.


\begin{thebibliography}{99}
\bibitem{BT82} E. Bedford and B. A. Taylor, A new capacity for plurisubharmonic functions, \emph{Acta Math.} \textbf{149} (1982), 1--40.
\bibitem{BT87} E. Bedford and B. A. Taylor, Fine topology, \v Silov boundary, and $(\ddc)^n$, \emph{J. Funct. Anal.} \textbf{72} (1987), 225--251.
\bibitem{BM22} F. Bianchi and S. Mongodi, On minimal kernels and Levi currents on weakly complete complex manifolds, \emph{Proc. Amer. Math. Soc.} \textbf{150} (2022), 3927--3939.
\bibitem{BM24} F. Bianchi and S. Mongodi, Generalized Levi currents and singular loci for families of plurisubharmonic functions, \emph{J. Geom. Anal.} \textbf{34} (2024), article 203.
\bibitem{BirkenhakeLange} C. Birkenhake and H. Lange, \emph{Complex Abelian Varieties}, 2nd ed., Springer, 2004.
\bibitem{Blocki} Z. B\l ocki, A note on maximal plurisubharmonic functions, \emph{Uzbek Math. J.} \textbf{1} (2009), 28--32.
\bibitem{Demailly82} J.-P. Demailly, Estimations $L^2$ pour l'op\'erateur $\bar\partial$ d'un fibr\'e vectoriel holomorphe semi-positif au-dessus d'une vari\'et\'e k\"ahl\'erienne compl\`ete, \emph{Ann. Sci. \'Ecole Norm. Sup. (4)} \textbf{15} (1982), 457--511.
\bibitem{DieuUnicity} N. Q. Dieu, A unicity theorem for plurisubharmonic functions, \emph{Ann. Polon. Math.} \textbf{100} (2011), 159--165.
\bibitem{DinewGuedjZeriahi} S. Dinew, V. Guedj and A. Zeriahi,
Open problems in pluripotential theory,
\emph{Complex Var. Elliptic Equ.} \textbf{61} (2016), no.~7, 902--930.
\bibitem{DoHiepComp} T. D. Do and P. H. Hiep,
Comparison principles for local plurisubharmonic potentials on complex manifolds and applications,
preprint (2026), arXiv:2610.05948.
\bibitem{DoHiep2026} T. D. Do and P. H. Hiep,
Local-to-global maximality for plurisubharmonic functions with locally analytic singularities,
preprint (2026), arXiv:2610.05937.
\bibitem{MST18} S. Mongodi, Z. S\l odkowski and G. Tomassini, Weakly complete complex surfaces, \emph{Indiana Univ. Math. J.} \textbf{67} (2018), 899--935.
\bibitem{Rashkovskii1998} A. Rashkovskii, Maximal plurisubharmonic functions associated with holomorphic mappings, \emph{Indiana Univ. Math. J.} \textbf{47} (1998), 297--310.
\bibitem{Shcherbina} N. V. Shcherbina, On compact subsets possessing strictly plurisubharmonic functions, \emph{Izv. Math.} \textbf{85} (2021), 605--618.
\bibitem{Siu1974} Y.-T. Siu, Analyticity of sets associated to Lelong numbers and the extension of closed positive currents, \emph{Invent. Math.} \textbf{27} (1974), 53--156.
\bibitem{Slodkowski} Z. S\l odkowski, On compact pseudoconcave sets, \emph{Proc. Amer. Math. Soc.} \textbf{152} (2024), 2449--2461.
\end{thebibliography}
\end{document}